\documentclass{amsart}
\pdfoutput=1
\usepackage{epsfig}
\usepackage[T1]{fontenc}    
\usepackage{ae,aecompl}     
\usepackage{amsmath}        
\usepackage{amssymb}
\usepackage{amsthm}
\usepackage{amsfonts}
\usepackage{amscd}
\usepackage{color}
\usepackage{mathrsfs}       
\usepackage{bbm}            
\usepackage{bm}             
\usepackage{url}            
\usepackage{paralist}       
\usepackage{graphics}       
\usepackage{epic}          
\usepackage[rel]{overpic}   
\usepackage{psfrag}
\usepackage{xspace}         

\newtheorem{thm}{Theorem}[section]
\newtheorem{cor}[thm]{Corollary}

\newtheorem{prop}[thm]{Proposition}
\newtheorem{rem}[thm]{Remark}

\newtheorem{expl}[thm]{Example}

\def\Ass{{\mathsf{Asso}}}  
\def\Perm{{\mathsf{Perm}}}
\def\H{{\mathscr H}}
\def\ID{{\operatorname {Id}}}
\def\a{{\alpha}}

\newcommand{\scalprod}[2]{{\langle #1,#2\rangle}}                         

\begin{document}


\title{Isometry classes of  generalized associahedra}

\author[N. Bergeron]{Nantel~Bergeron}
\address[Nantel Bergeron]{Department of Mathematics and Statistics\\
York University\\
Toronto, Ontario, M3J 1P3\\
CANADA}
\email{bergeron@mathstat.yorku.ca}

\author[C. Hohlweg]{Christophe~Hohlweg}
\address[Christophe Hohlweg]{Universit\'e du Qu\'ebec \`a Montr\'eal\\
LaCIM et D\'epartement de Math\'ematiques\\ CP 8888 Succ. Centre-Ville\\
Montr\'eal, Qu\'ebec, H3C 3P8\\ CANADA}
\email{hohlweg.christophe@uqam.ca}
\urladdr{http://www.lacim.uqam.ca/\~{}hohlweg}

\author[C. Lange]{Carsten~Lange}
\address[Carsten Lange]{Freie Universit{\"a}t Berlin\\
Fachbereich f{\"u}r Mathematik und Informatik\\
Arnimallee 3\\
14195 Berlin\\
GERMANY}
\email{lange@math.tu-berlin.de}

\author[H. Thomas]{Hugh~Thomas$^*$}
\address[Hugh Thomas]{Department of Mathematics and Statistics\\
University of New Brunswick\\
Fredericton, New Brunswick, E3B 5A3\\
CANADA}
\email{hugh@math.unb.ca}

\date{April 16, 2008}

\thanks{$^*$ This work is supported in part by CRC, NSERC and PAFARC.
It is the results of a working session at the Algebraic
Combinatorics Seminar at the Fields Institute with the active
participation of Anouk Bergeron-Brlek, Philippe Choquette, Huilan
Li, Trueman Machenry, Alejandra Premat, Muge Taskin and Mike Zabrocki.}

\begin{abstract}
\noindent Let $(W,S)$ be a finite Coxeter system acting by
reflections on an $\mathbb R$-Euclidean space with simple roots
$\Delta=\{\a_s\,|\, s\in S\}$ of the same length,
 and fundamental weights $\Delta^*=\{v_s\,|\,
s\in S\}$. We set $M(e)=\sum_{s\in S}\kappa_s v_s$, $\kappa_s>0$
and for $w\in W$ we set $M(w)=w\big(M(e)\big)$. The permutahedron
$\Perm(W)$ is the convex hull of the set $\{M(w)\,|\, w\in W\}$.
Given a Coxeter element $c\in W$, we have defined in a previous
work a generalized associahedron $\Ass_c(W)$ whose normal fan is
the corresponding $c$-Cambrian fan $\mathcal F_c$ defined by
N.~Reading. By construction, $\Ass_c(W)$ is obtained  from
$\Perm(W)$ by removing some halfspaces according to a rule
prescribed by $c$. In this work, we classify  the isometry classes
of these realizations. More precisely, for $(W,S)$ an irreducible
finite Coxeter system and $c,c'$ two Coxeter elements in $W$, we
have that $\Ass_{c}(W)$ and $\Ass_{c'}(W)$ are isometric if and
only if $\mu(c') = c$ or~$\mu(c')=w_0 c^{-1}w_0$ for $\mu$ an
automorphism of the Coxeter graph of $W$ such that
$\kappa_s=\kappa_{\mu(s)}$ for all $s\in S$. As a byproduct, we
classify the isometric Cambrian fans of $W$.
 \end{abstract}

 \maketitle


\section{Introduction.}\label{se:Intro}

Studying binary operations, Stasheff \cite{stasheff, stasheff2}
remarked that there was a cell complex whose vertices correspond
to the possible compositions of $n$ binary operations.
Furthermore this cell
complex
can be realized as a simple polytope, the associahedron. There is a
natural relation between the permutahedron (weak lattice on
permutations) and the associahedron: the permutahedron can
naturally be written as an intersection of halfspaces indexed by
weights for the $A_n$ root system. If one intersects a certain,
carefully chosen subset of these halfspaces, one can obtain the
associahedron (see Example~\ref{permass} below).
Shnider-Sternberg~\cite{shnider_sternberg} and Loday~\cite{loday}
 give us a beautiful explicit construction of the associahedron
from the permutahedron. The facets of both objects encode the
algebraic structure of the corresponding free operads.

Generalized associahedra were introduced by S.~Fomin and
A.~Zelevinsky in their work on cluster algebras
\cite{fomin_zelevinsky}. The geometry of these objects encode nice
algebraic structures. Therefore one important question is to find
good polytopal realizations of the generalized associahedra. This
was first answered in \cite{chapoton_fomin_zelevinsky}. Then,
N.~Reading \cite{reading4} constructed a family of fans, the
Cambrian fans $\{\mathcal F_c\}$ indexed by Coxeter elements $c$
of a given finite Coxeter group $W$. More recently, we have
constructed \cite{realisation2} a family of generalized
associahedra $\Ass_{c}(W)$, one for each Cambrian fan $\mathcal
F_c$. These generalized associahedra are realized from the
corresponding
permutahedron by removing some halfspaces according to a rule
specified by $c$, linking questions about generalized associahedra
to questions about the better known permutahedron.

It is now natural to ask how many distinct (up to isometry)
generalized associahedra we have. This is what we answer here.
Our main theorem (Theorem~\ref{thm:Main}) is to describe
completely the isometry classes of generalized associahedra as
realized in \cite{realisation2} . The isometry classes depend of
the choice of the starting permutahedron. As a byproduct we obtain
a classification of the isometry classes of Cambrian fans
(Corollary~\ref{cor:Cambrian}):  the Cambrian fans indexed by
Coxeter elements $c$ and $c'$ are isometric if and only if
$\mu(c^{\prime})=c$ or $\mu(c^{\prime})=c^{-1}$ for some $\mu$ an
automorphism of the Coxeter graph of $W$. In Section~2 we
introduce the necessary definitions and state our main theorem
(Theorem~\ref{thm:Main}). The proof is found in Section~4.
Section~3 is dedicated to some auxiliary results needed for this
proof.  For most of the paper, we make the simplifying assumption
that the Coxeter system in question is irreducible; in Section~5, we
explain how to deal with the reducible case.

\section{background and  main theorem.}\label{se:main}
We assume some basic familiarity with Coxeter groups and root
systems and follow the notation of~\cite{humphreys}. Let~$(W,S)$
be a finite Coxeter system acting by reflections on an $\mathbb
R$-euclidean space~$(V,\scalprod{\cdot}{\cdot})$ with length
function $\ell: W\rightarrow \mathbb N$. Without loss of
generality, we assume that the action of $W$ is essential relative
to $V$, that is,  has no non trivial space fixed pointwise.

 Let $\Phi$ be  a root system corresponding to $(W,S)$, with all roots having
equal length.  (In particular, we do not assume $\Phi$ is crystallographic.)
The simple roots $\Delta$ are a basis of $V$, and the reflection $s$ maps
$\a_s$ to $-\a_s$ and fixes the
 hyperplane $H_s=\{v\in V\,|\, \scalprod{v}{\a_s}=0\}$.
 Let $\Delta^*=\{v_s\,|\, s\in S\}$ be the set of fundamental weights of $\Delta$, that is, $\scalprod{v_t}{\a_s}=1$ if $s=t$ and $0$ otherwise.
 As $V$ is finite dimensional we identify $V$ and $V^*$.

We now aim for a definition of the $W$-permutahedron and pick a point $u \in V$ contained in the complement of the
reflection hyperplanes of~$W$. Without loss of generality, we choose
\[
  u := \sum_{s\in S} \kappa_s v_s, \qquad \kappa_{s} > 0.
\]
For $w \in W$ we write
\[
  M(e) := u \qquad\text{and}\qquad M(w) := w\left( M(e) \right)
\]
and obtain the permutahedron $\Perm_u(W)$ as convex hull of $\{M(w)\,|\, w\in W\}$. The index~$u$ will often be
omitted for brevity. Equivalently, we have
\[
  \Perm(W) = \bigcap_{s\in S}\
             \bigcap_{x\in W}\H_{(x,s)}
\]
where
\[
   \H_{(x,s)}
    := \{ v \in V\,|\,
           \scalprod{v}{x(v_{s})}\leq \scalprod{M(e)}{v_{s}}\}.
\]
We also make use of the hyperplane~$H_{(x, s)}=\{ v \in V\,|\, \scalprod{v}{x(v_{s})}= \scalprod{M(e)}{v_{s}}\}$.

Denote by $W_I$ the {\em standard parabolic subgroup of $W$}
generated  by $I\subseteq S$. Note that $H_{(w,s)}=H_{(x,s)}$ if
and only if $w\in xW_{S\setminus \{s\}}$. Also, $M(w)\in
H_{(x,s)}$ if and only if $H_{(x,s)}=H_{(w,s)}$. Hence we have a
simple way to describe the vertices:
\[
   \{M(w) \} = \bigcap_{s\in S} H_{(w,s)} .
\]

\begin{figure}
  \begin{center}
  \begin{minipage}{0.95\linewidth}
     \begin{center}
     \begin{overpic}
        [width=140pt]{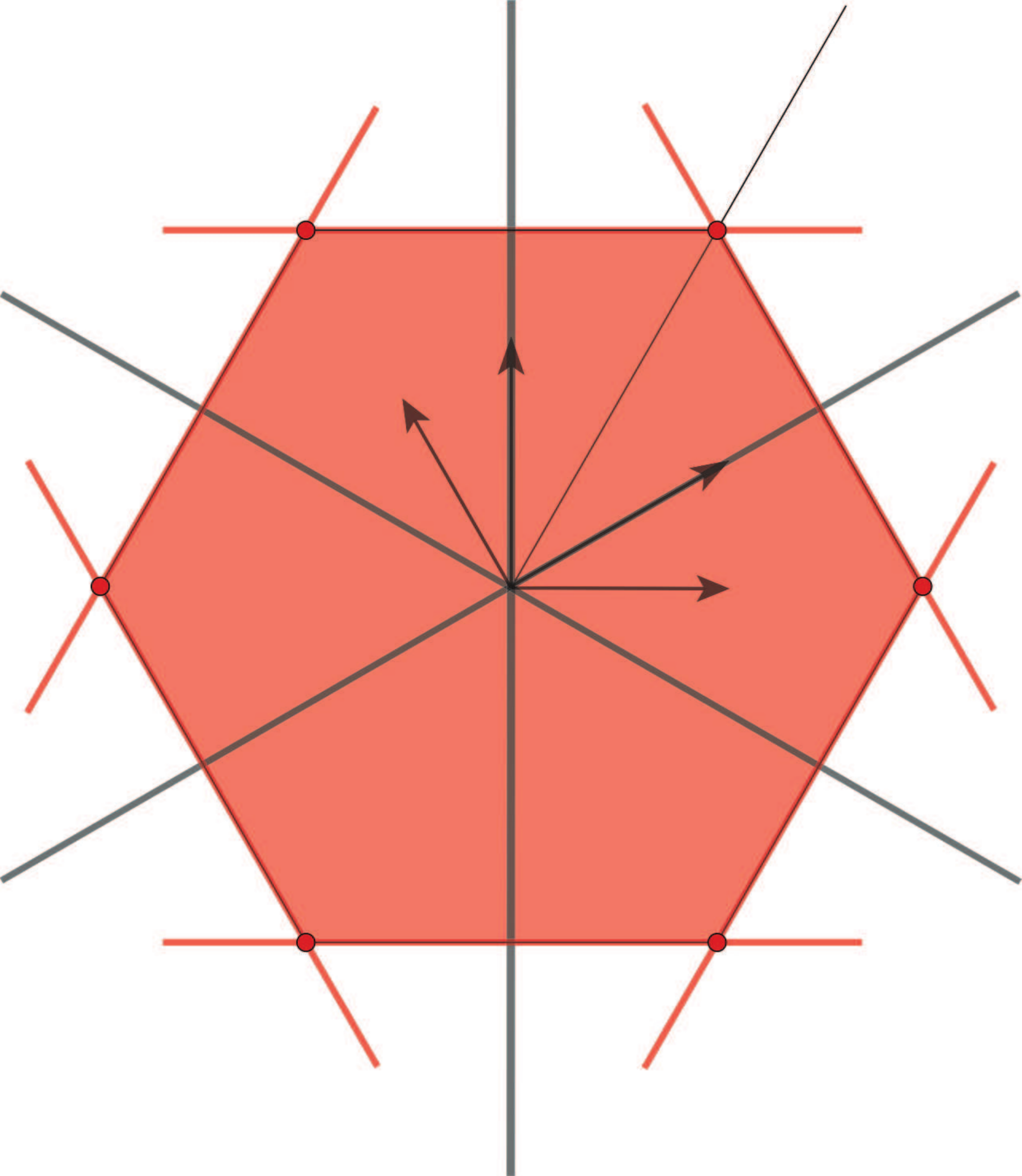}
        \put(71,95){\tiny $L$}
        \put(63,82){\tiny $M(e)$}
        \put(44,95){{\tiny {$H_{s_1}$}}}
        \put(13,82){\tiny $M(s_1)$}
        \put(-10,49){\tiny $M(s_1s_2)$}
        \put(81,49){\tiny $M(s_2)$}
        \put(62,16){\tiny $M(s_2s_1)$}
        \put(5,15){\tiny $M(s_1s_2s_1)$}
        \put(80,68){\tiny $H_{s_2}$}
        \put(-3,76){\tiny $H_{s_1s_2s_1}$}
        \put(29,65){\tiny $\a_2$}
        \put(62,48){\tiny $\a_1$}
        \put(63,59){\tiny $v_{s_1}$}
        \put(45,68){\tiny $v_{s_2}$}
        \put(85,39){\tiny $H_{(e,s_1)}=H_{(s_2,s_1)}$}
        \put(85,59){\tiny $H_{(s_2,s_2)}=H_{(s_2s_1,s_2)}$}
        \put(74,79){\tiny $H_{(e,s_2)}=H_{(s_1,s_2)}$}
        \put(-44,39){\tiny $H_{(s_1,s_1)}=H_{(s_1s_2,s_1)}$}
        \put(-53,59){\tiny $H_{(s_1s_2,s_2)}=H_{(s_1s_2s_1,s_2)}$}
        \put(-41,20.5){\tiny $H_{(s_1s_2s_1,s_1)}=H_{(s_2s_1,s_1)}$}
     \end{overpic}
     \end{center}
     \caption[]{The permutahedron~$\Perm(S_3)$ obtained as convex hull of the $S_3$-orbit
                of~$M(e) \in L$ or as intersection of the half spaces~$\H_{(x,s)}$.}
     \label{fig:construction_perm}
  \end{minipage}
  \end{center}
\end{figure}

\begin{expl}[Realization of~$\Perm(A_2)$]\textnormal{
   We consider the Coxeter group~$W=S_3$ of type~$A_2$ acting on~$\mathbb R^2$. The reflections $s_1$ and $s_2$
   generate~$W$. The simple roots that correspond to~$s_1$ and~$s_2$ are~$\a_1$ and~$\a_2$. They are normal to
   the reflection hyperplanes~$H_{s_1}$ and~$H_{s_2}$. The dual vectors to the simple roots correspond to the vectors~$v_1$
   and~$v_2$. Fix a ray~$L = \{ \mu (\kappa_1 v_1 + \kappa_2 v_2) \;|\; \mu >0\}$ where
   $\kappa_1,\kappa_2 > 0$. We choose $M(e) \in L$ and obtain the permutahedron as convex hull of the $W$-orbit
   of~$M(e)$. Alternatively, the permutahedron can be described as intersection of the half spaces~$\H_{(x,s)}$
   with bounding hyperplanes~$H_{(x,s)}$ for $x\in W$ and $s \in S$. All the objects are indicated in
   Figure~\ref{fig:construction_perm}.}
\end{expl}

\smallskip

For $c$ a Coxeter element in $W$,  that is to say, the product of
the simple reflections of $W$ taken in some order, and $I\subseteq
S$, we denote by $c_{(I)}$ the subword of $c$ obtained by taking
only the simple reflections in $I$. So $c_{(I)}$ is a Coxeter element
of $W_I$. Reading defined the {\em $c$-sorting word} of $w\in W$
in
 \cite[Section~$2$]{reading2} as the unique subword of the
infinite word $c^\infty=cccccc\dots$  that is a reduced expression
for~$w$ and is the lexicographically smallest sequence of
positions occupied by this subword. In particular, the $c$-sorting
word of~$w$ is such that $w=c_{(K_1)}c_{(K_2)}\dots c_{(K_p)}$ with
non-empty~$K_i \subseteq S$ and  $\ell(w)=\sum_{i=1}^p |K_i|$. As
example we consider the Coxeter group~$W=S_4$ of type~$A_3$
generated by the simple reflections~$S = \{ s_1, s_2, s_3 \}$,
where $s_1,s_3$ commute, and the Coxeter element~$c=s_2s_1s_3$.
The $c$-sorting word of the longest element~$w_0\in W$ is
$s_2s_1s_3s_2s_1s_3=c_{(S)}c_{(S)}$. If we choose the Coxeter
element~$c=s_1s_2s_3$ instead of~$s_2s_1s_3$, then the $c$-sorting
word of~$w_0$
is~$s_1s_2s_3s_1s_2s_1=c_{(S)}c_{(\{s_1,s_2\})}c_{(\{s_1\})}$.

The sequence~$c_{(K_1)}, \ldots,c_{(K_p)}$ associated to the
$c$-sorting word for~$w$ is called {\em the $c$-factorization}
of~$w$. The $c$-factorization of~$w$ is independent of the chosen
reduced word for~$c$ but depends on the Coxeter element~$c$. In
general the $c$-factorization does not yield a nested sequence $K_1,
\ldots, K_p$ of subsets of~$S$. An element~$w \in W$ is called {\em
$c$-sortable} if $K_1 \supseteq K_2 \supseteq \ldots \supseteq K_p$.
Reading proves in~\cite{reading2} that the longest element $w_0 \in
W$ is $c$-sortable for any chosen Coxeter element~$c$.

Given a specific reduced
word $\bf v$, we say that $u$ is a prefix up to commutation
of $\bf v$ if some reduced word for $u$
appears as a prefix of a word which can be obtained
from $\bf v$ by the commutation of commuting reflections.
In~\cite{realisation2} we define an element~$w \in W$
to be a {\em $c$-singleton} if it is a prefix up to commutation
of the $c$-factorization
of $w_0$.
We illustrate this notion by
considering again the Coxeter group~$W=S_4$ and the Coxeter
element $c=s_2s_3s_1$. The $c$-singletons are
\begin{center}
\begin{tabular}{lll}
   $e$,            & $s_2s_3$,             & $s_2s_1s_3s_2s_1$,\\
   $s_2$,          & $s_2s_1s_3$,          & $s_2s_1s_3s_2s_3$, and\\
   $s_2s_1,\qquad$ & $s_2s_1s_3s_2,\qquad$ & $w_0=s_2s_3s_1s_2s_3s_1$.
\end{tabular}
\end{center}
For example $s_2s_1$ is a not a prefix of the $c$-factorization
$w_0=
s_2s_3s_1s_2s_3s_1$,
but it is a prefix up to commutation because it appears as a prefix
after commuting the simple reflections $s_1$, $s_3$.

The halfspace~$\H_{(x,s)}$ is said to be {\em $c$-admissible} if
the hyperplane~$H_{(x,s)}$ contains~$M(w)$ for some
$c$-singleton~$w$. We have shown in~\cite{realisation2} that the
intersection of all $c$-admissible halfspaces~$\H_{(x,s)}$ is a
{\em generalized associahedron}~$\Ass_{c}(W)$ whose normal fan is
the $c$-Cambrian fan $\mathcal F_c$ (see \cite{reading4} for a
definition of~$\mathcal F_c$).

\begin{expl}\label{permass}\textnormal{
   The Coxeter group~$W=S_3$ generated by the reflections~$s_1$ and~$s_2$ has two Coxeter
   elements: $c_1=s_1s_2$ and $c_2=s_2s_1$. The $c_1$-singletons are $e$, $s_1$, $s_1s_2$, and
   $s_1s_2s_1$ while the  $c_2$-singletons are $e$, $s_2$, $s_2s_1$,and $s_2s_1s_2$. Starting
   with the permutahedron~$\Perm(S_3)$, we obtain the two associahedra~$\Ass_{c_1}(S_3)$ and
   $\Ass_{c_2}(S_3)$ shown in Figure~\ref{fig:A_2_associahedra} as intersection of the $c_1$-
   and $c_2$-admissible halfspaces.}
\end{expl}
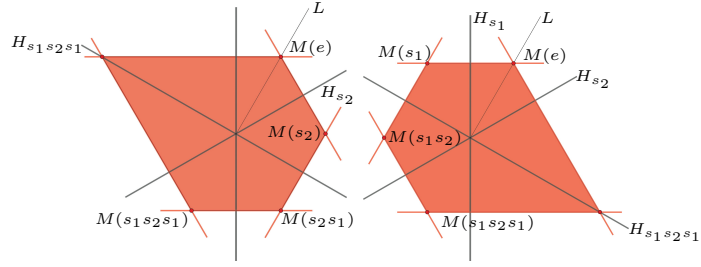
\begin{figure}
  \begin{center}
  \begin{minipage}{0.95\linewidth}
     \begin{center}
     \begin{overpic}
        [width=100pt]{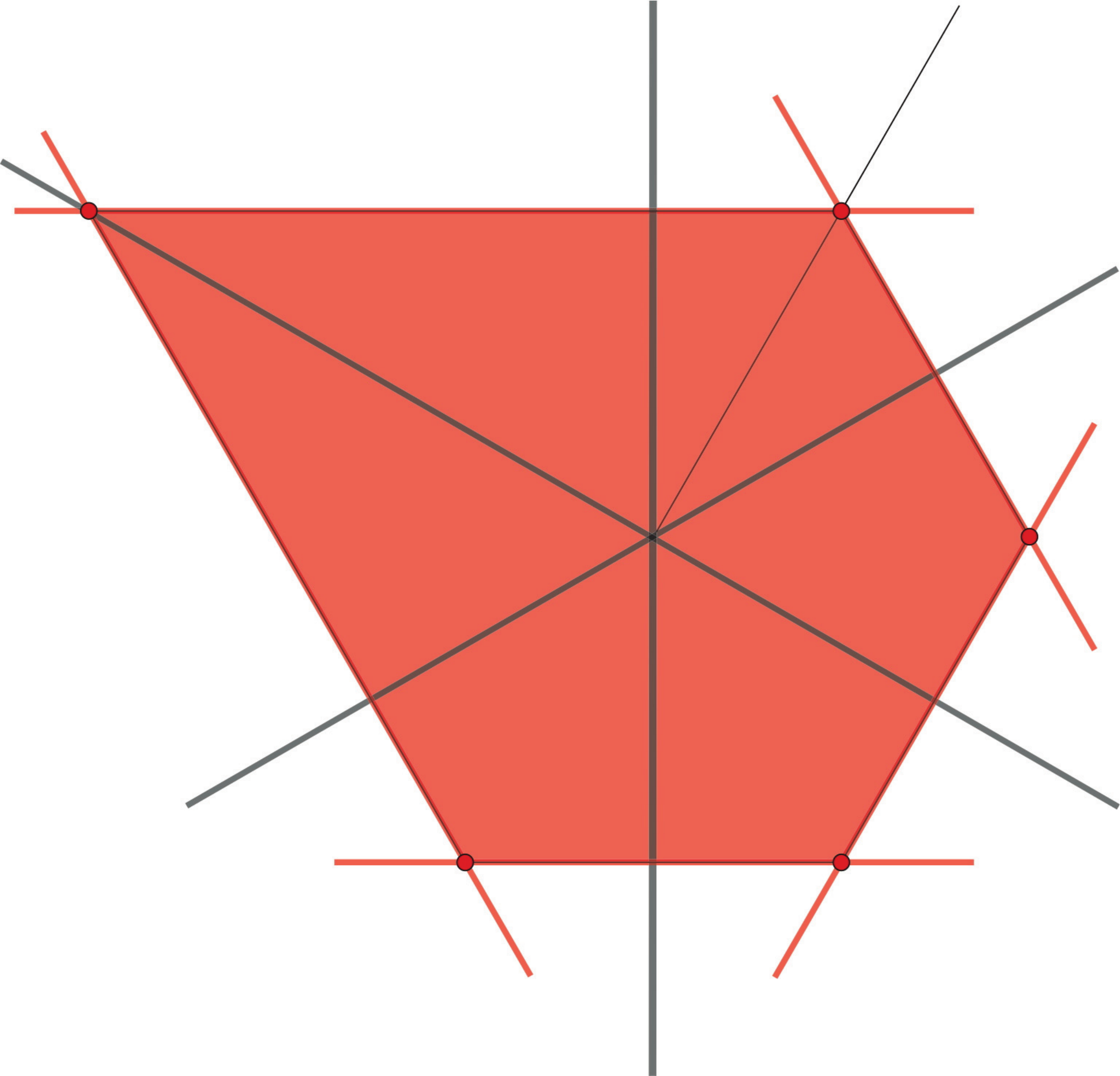}
        \put(87,94){\tiny $L$}
        \put(78,79){\tiny $M(e)$}
        \put(70,47){\tiny $M(s_2)$}
        \put(75,14){\tiny $M(s_2s_1)$}
        \put(5,14){\tiny $M(s_1s_2s_1)$}
        \put(90,62){\tiny $H_{s_2}$}
        \put(-27,81){\tiny $H_{s_1s_2s_1}$}
     \end{overpic}\ \
     \begin{overpic}
        [width=100pt]{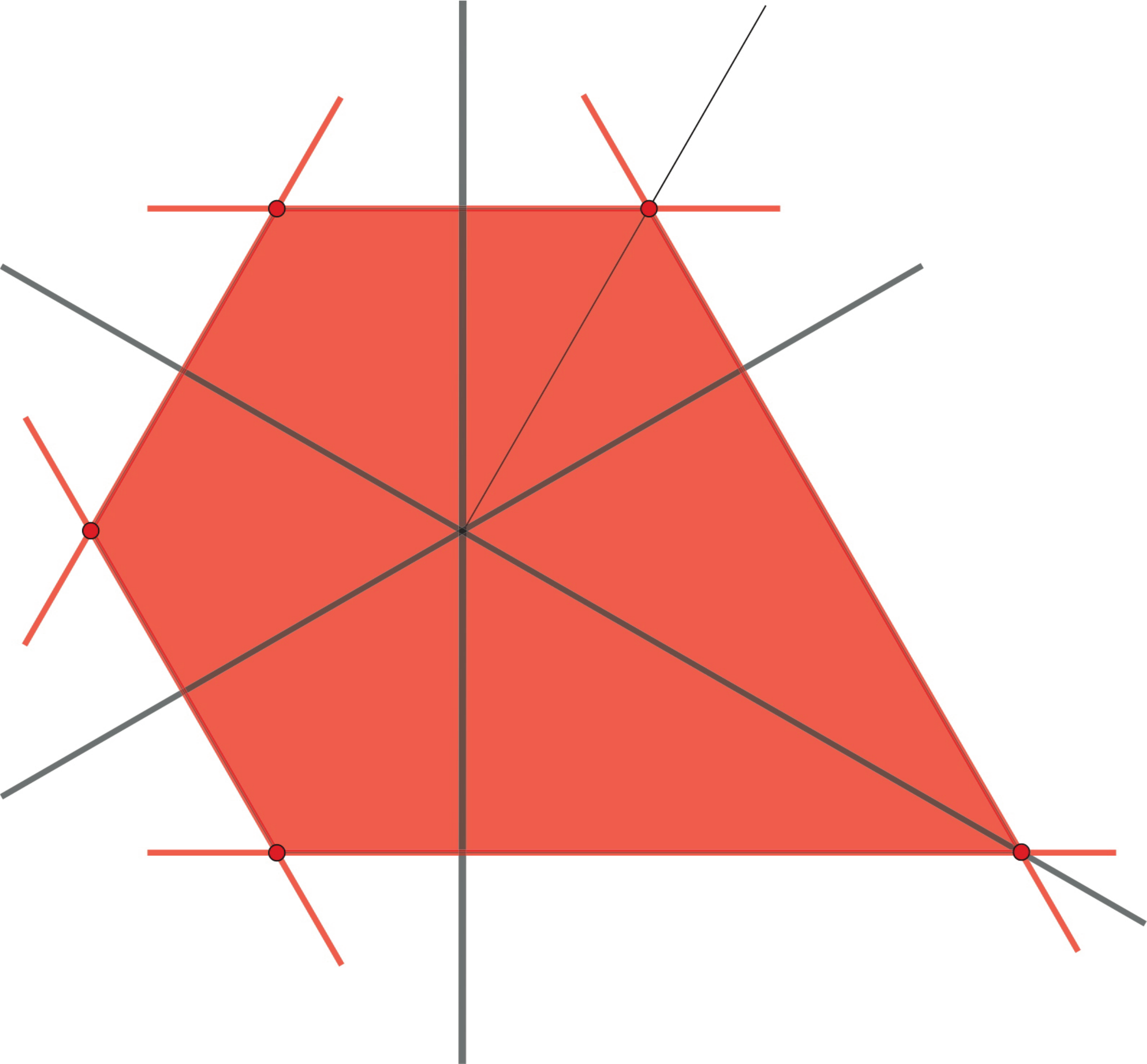}
        \put(67,90){\tiny $L$}
        \put(60,76){\tiny $M(e)$}
        \put(41,90){{\tiny {$H_{s_1}$}}}
        \put(5,77){\tiny $M(s_1)$}
        \put(9,45){\tiny $M(s_1s_2)$}
        \put(29,13){\tiny $M(s_1s_2s_1)$}
        \put(80,68){\tiny $H_{s_2}$}
        \put(100,10){\tiny $H_{s_1s_2s_1}$}
     \end{overpic}
     \end{center}
     \caption[]{The two associahedra~$\Ass_{c_1}(S_3)$ (left) and~$\Ass_{c_2}(S_3)$ (right) obtained from the
                permutahedron~$\Perm(S_3)$ by keeping the $c$-admissible halfspaces~$\H_{(x,s)}$.}
     \label{fig:A_2_associahedra}
  \end{minipage}
  \end{center}
\end{figure}

For most of the paper, we will assume that $(W,S)$ is irreducible.
The case where $(W,S)$ is reducible requires a straightforward (but not
immediate) extension of the results in the
irreducible case; we describe this in the final
section.

An automorphism of the Coxeter graph associated to~$(W,S)$ is a
bijection~$\mu$ on~$S$ such that the order of~$\mu(s)\mu(t)$
equals the order of $st$ for all $s,t \in S$. In particular, $\mu$
induces an automorphism on~$W$.

Let  $u=\sum_{s\in S}\kappa_s v_s$ be a point in $V$. We will say
that {\em $u$ is balanced} if $\kappa_s=\kappa_t$ for all $s,t\in
S$. An automorphism $\mu$ of the Coxeter graph is a {\em
$u$-automorphism} if $\kappa_s=\kappa_{\mu(s)}$ for all $s\in S$.
 In particular, if $u$ is balanced, then any automorphism of a Coxeter graph is a $u$-automorphism.

\begin{thm}\label{thm:Main} Let~$(W,S)$ be an irreducible finite Coxeter system and $c_1$, $c_2$ be two Coxeter elements
 in $W$. Suppose that~$u=\sum_{s \in S} \kappa_s v_s$ for some~$\kappa_s > 0$. 
  The following statements are equivalent.
   \begin{compactenum}
      \item $\Ass_{c_1}(W)=\varphi\left(\Ass_{c_2}(W)\right)$ for some linear isometry~$\varphi$ on~$V$.
      \item There is an $u$-automorphism~$\mu$ of the Coxeter graph of~$(W,S)$ such that
       $\mu(c_2) = c_1$~or~$\mu(c_2)=w_0c_1^{-1}w_0$.
   \end{compactenum}
\end{thm}

Observe that $w_0c^{-1}w_0$ may or may not equal $c$ (for instance in $A_3$ take $c=s_1s_3s_2$).
 So the second condition in Theorem~\ref{thm:Main} may be
redundant and the associahedra may actually be identical (not just isometric).
Moreover, if the $\kappa_s$ are chosen generically, that is
distinct, then the isometry classes are of cardinality $1$ or $2$.
As stated in the next corollary, the isometry classes reach their
maximal cardinality if $u$ is balanced.

\begin{cor}\label{cor:Main} Let~$(W,S)$ be an irreducible finite Coxeter system and $c_1$, $c_2$ be two Coxeter elements
 in $W$. If $u$ is balanced, then the following statements are equivalent.
   \begin{compactenum}
      \item $\Ass_{c_1}(W)=\varphi\left(\Ass_{c_2}(W)\right)$ for some linear isometry~$\varphi$ on~$V$.
      \item There is an automorphism~$\mu$ of the Coxeter graph of~$(W,S)$ such that $\mu(c_2) = c_1$ or~$\mu(c_2)=c_1^{-1}$.
   \end{compactenum}
\end{cor}
\begin{proof} It follows from Theorem~\ref{thm:Main} and a
rewriting of the second assumption in accordance with the fact that
the map $s\mapsto w_0 s w_0$ is an automorphism of the Coxeter
graph.
\end{proof}

If $u$ is balanced, then
$\Ass_{c_1}(W)=\varphi\left(\Ass_{c_2}(W)\right)$ for some linear
isometry~$\varphi$ on~$V$ if and only if there is an
$u$-automorphism~$\theta$ of the Bruhat ordering of~$(W,S)$ such
that~$\theta(c_2)=c_1$. This follows by inspection for
$|S|=1,2$ and, for $|S|\geq 3$, from a characterization of
automorphisms of Bruhat orderings due to van den Hombergh, see
Section~$8.8$ of~\cite{humphreys} and the fact that a Coxeter
element~$c$ defines an orientation of the Coxeter graph~$\Gamma$:
Orient the edge $\{s_i,s_j\}$ from~$s_i$ to~$s_j$ if and only
if~$s_i$ is to the left of~$s_j$ for any reduced word for~$c$.
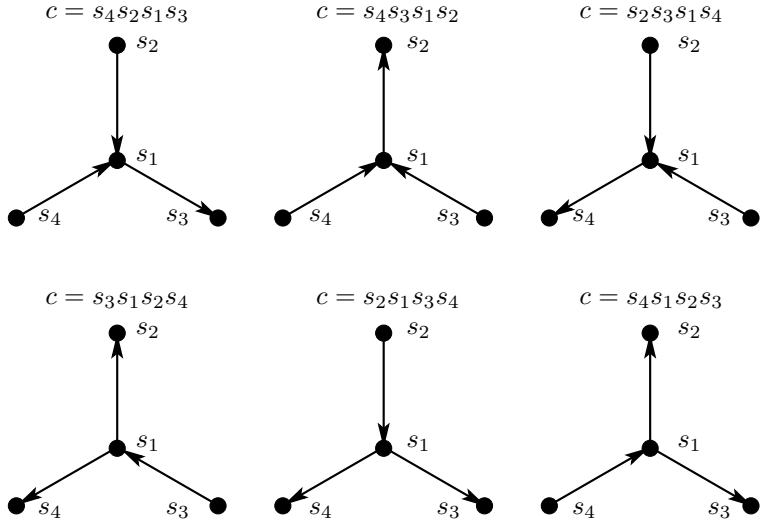
\begin{figure}[h]
  \begin{center}
  \begin{minipage}{0.95\linewidth}
     \begin{center}
     \begin{overpic}
        [width=10cm]{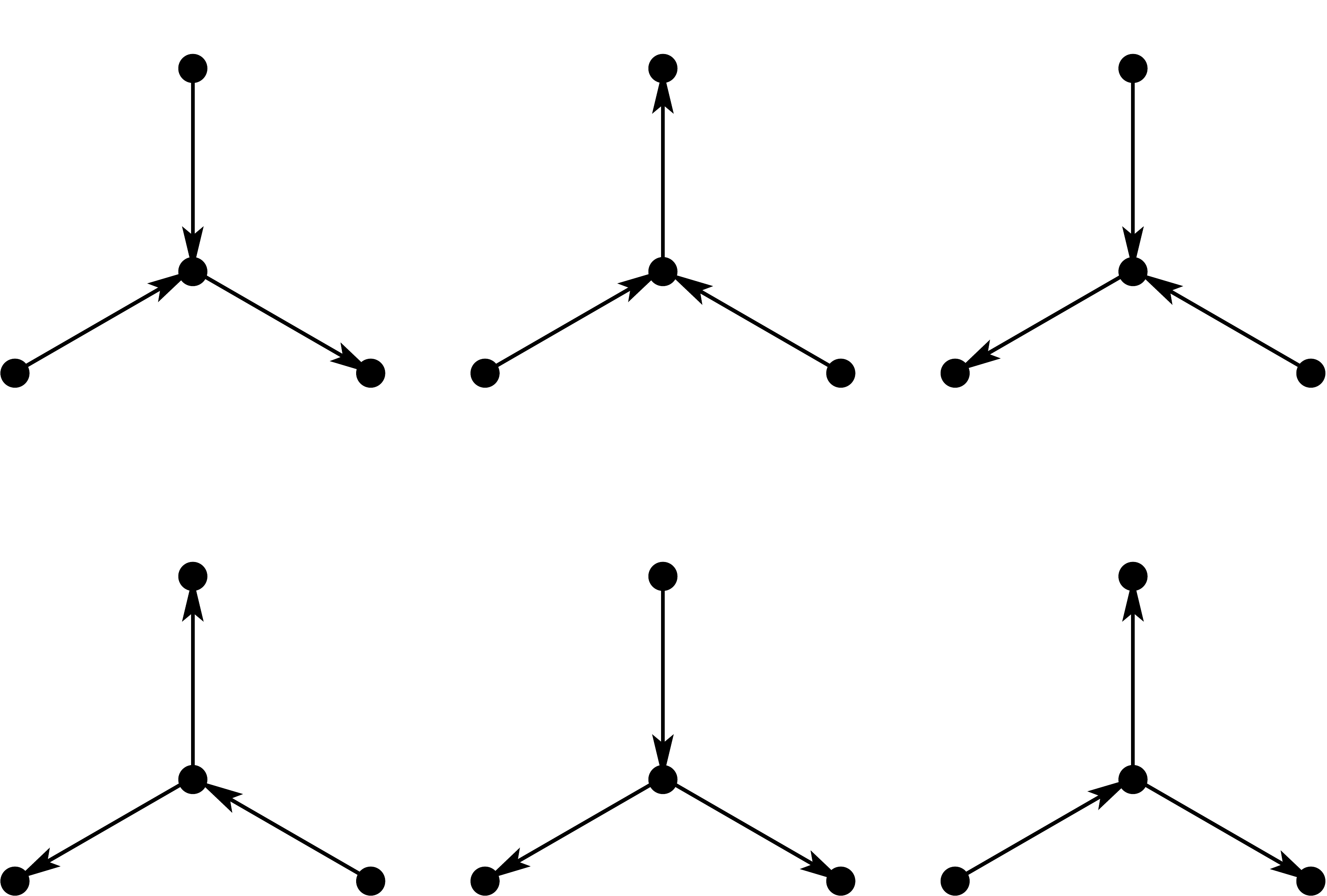}
        \put(17,47){$s_1$}
        \put(17,62){$s_2$}
        \put(21,39){$s_3$}
        \put(4,39){$s_4$}
        \put(5,66){$c=s_4s_2s_1s_3$}
        \put(53,47){$s_1$}
        \put(53,62){$s_2$}
        \put(57,39){$s_3$}
        \put(40,39){$s_4$}
        \put(41,66){$c=s_4s_3s_1s_2$}
        \put(89,47){$s_1$}
        \put(89,62){$s_2$}
        \put(93,39){$s_3$}
        \put(75,39){$s_4$}
        \put(76,66){$c=s_2s_3s_1s_4$}
        \put(17,9){$s_1$}
        \put(17,24){$s_2$}
        \put(21,0){$s_3$}
        \put(4,0){$s_4$}
        \put(5,28){$c=s_3s_1s_2s_4$}
        \put(53,9){$s_1$}
        \put(53,24){$s_2$}
        \put(57,0){$s_3$}
        \put(40,0){$s_4$}
        \put(41,28){$c=s_2s_1s_3s_4$}
        \put(89,9){$s_1$}
        \put(89,24){$s_2$}
        \put(93,0){$s_3$}
        \put(75,0){$s_4$}
        \put(76,28){$c=s_4s_1s_2s_3$}
     \end{overpic}
     \end{center}
     \caption[]{Six Coxeter elements and their associated oriented Coxeter graphs of the Coxeter group of type~$D_4$
                that yield isometric associahedra.}
     \label{fig:six_orientations}
  \end{minipage}
  \end{center}
\end{figure}

\medskip
Theorem~\ref{thm:Main} combined with the classification of
irreducible finite Coxeter groups yields that the cardinality of
an isometry class in the case where $u$ is balanced is either two,
four, or six. We briefly discuss the situation.
\begin{expl}\textnormal{Use the notation and hypothesis of Corollary~\ref{cor:Main}.
   \begin{compactenum}
      \item Let $(W,S)$ be a Coxeter system of type~$A_n\ (n\geq 2)$, $E_6$, $F_4$, or~$I_2(m)$. Then there is precisely one non-trivial
            automorphism~$\mu$ of the Coxeter graph. Hence there are either two or four elements
            in the isometry class of~$\Ass_{c}(W)$. The cardinality equals two if $\mu(c)\in\{c,c^{-1}\}$ and
            equals four if $\mu(c)\not\in\{c,c^{-1}\}$.
      \item Let $(W,S)$ be a Coxeter system of type~$B$, $E_7$, $E_8$, $H_3$, or~$H_4$. Then the conjugation by $w_0$
            is the identity, and~$\ID$
            is the only automorphism of the associated Coxeter graph. So each isometry class has cardinality two,
            only the Coxeter elements~$c$ and~$c^{-1}$ yield isometric associahedra.
      \item Let~$(W,S)$ be a Coxeter system of type $D$. If $|S|>4$ then there is only one non-trivial automorphism~$\mu$
            of the Coxeter graph and the isometry class of~$\Ass_{c}(W)$ has cardinality two if $\mu(c)\in\{c,c^{-1}\}$
            and four otherwise.
            If $|S|=4$, the group of automorphisms of the Coxeter graph is generated by the
            non-trivial automorphisms~$\mu$ and~$\nu$ with $\mu^2=\ID$ and~$\nu^3=\ID$. The isometry class
            of~$\Ass_{c}(W)$ consists either of two or six elements, see Figure~\ref{fig:six_orientations}
            for six Coxeter elements that yield isometric associahedra.
   \end{compactenum}}
\end{expl}

   Theorem~\ref{thm:Main} allows a classification of the isometric Cambrian fans as
  well. The proof will be given at the end of Section~\ref{se:Proof}.

  \begin{cor}\label{cor:Cambrian} The following propositions are equivalent:
  \begin{compactenum}
  \item The Cambrian fans~$\mathcal F_c$ and~$\mathcal F_{c'}$ are isometric;
  \item $\Ass_{c}(W)$ and~$\Ass_{c'}(W)$ are isometric if $u$ is balanced;
  \item there is an automorphism~$\mu$ of the Coxeter graph of~$(W,S)$ such that $\mu(c') = c$~or~$\mu(c')=c^{-1}$.
 \end{compactenum}
\end{cor}

\begin{rem}\textnormal{ In fact, the condition that $u$ be balanced in
(2) above can be weakened to require only that $u$ satisfies
$\kappa_s=\kappa_{w_0sw_0}$, for all $s\in S$ (see
Eq.~\ref{eq:proofvalid} in the proof)}.\end{rem}

\section{Preliminary results}\label{se:prelim}
\noindent From now on, we fix $u= \sum_{s \in S}\kappa_s v_s$ for
some constants~$\kappa_s > 0$. Remember that $M(e)=u$. 

\medskip
\noindent
Let~$(W,S)$ be a Coxeter system and consider the polyhedron~$P$ defined by:
\[
  P:=\bigcap_{s\in S}\H_{(e,s)} \cap \bigcap_{s\in S} \H_{(w_0,s)}.
\]

\begin{rem}\textnormal{In fact,~$P$ is a full-dimensional convex polytope because~$W$
acts essentially on~$V$ and the cones $\bigcap_{s\in S}\H_{(e,s)}$
and $\bigcap_{s\in S} \H_{(w_0,s)}$ are strictly convex, pointed
with apex~$M(e)$ and~$M(w_0)$, and both contain~$\Perm(W)$. In
other words,~$P$ is obtained from~$\Perm(W)$ by removing, from the
definition of $\Perm(W)$ as an intersection of halfspaces, all
halfspaces~$\H_{(x,s)}$ that satisfy $M(e) \not\in H_{(x,s)}$ and
$M(w_0) \not\in H_{(x,s)}$.}
\end{rem}

\begin{prop}\label{prop:UnicityPerm}
   Let $\varphi:V\to V$ be a linear isometry that maps~$P$ to itself and has the fixed point~$M(e)$.
   Then $\varphi$ induces an $u$-automorphism~$\mu$ of the Coxeter graph of~$(W,S)$ such that
   $\varphi(v_s):=v_{\mu(s)}$ for every $s \in S$.
\end{prop}

\begin{proof}
  For $\varphi$ satisfying our hypothesis, we have that $\varphi(M(e))=M(e)\in H_{(e,s)}$.
Hence $\varphi$ induces a bijection on the set $\{\H_{(e,s)}\,|\,
s\in S\}$. Since $v_s$ is a normal vector to $\H_{(e,s)}$ for any
$s\in S$, we have that $\varphi(v_s)=k_sv_{t_s}$ for some $k_s>0$
and $t_s\in S$. Hence $\sum_{s\in S}k_s\kappa_{s}
v_{t_s}=\varphi(M(e))=M(e)=\sum_{s\in S} \kappa_s v_s$ and then
$\kappa_{t_s}=k_s\kappa_s$.

On the other hand, since $\varphi$ is an isometry fixing $P$ and $M(e)$, it induces a bijection on the set of
 edges of $P$ which have $M(e)$ as one of their vertices. Each of these edges is contained in a
 line $l_s:=\bigcap_{r\in S\setminus\{s\}} H_{(e,r)}$, for $s\in S$. So $\varphi$ induces a bijection on the
 set $\{l_s\,|\, s\in S\}$.   Since $v_r$ is a fundamental weight of $\Delta$, and is a normal vector
 of $H_{(e,r)}$,  the hyperplane $H_{(e,r)}$ is spanned by the simple roots   $\{\alpha_u\,|\,u\in S\setminus\{r\}\}$,
 and thus $l_s$ is spanned by the simple root $\alpha_s$. As the simple roots are all
  of the same length and $\varphi$ preserves the norm of vectors, $\varphi(\alpha_s)=\pm  \alpha_{r_s}$ for some $r_s\in S$.
 Now
$$
1=\scalprod{v_s}{\alpha_s}=\scalprod{\varphi(v_s)}{\varphi(\alpha_s)}=\pm k_s\scalprod{v_{t_s}}{\alpha_{r_s}}.
$$

We conclude that $r_s=t_s$, $k_s=1$ and $\kappa_s=\kappa_{t_s}$ for all $s\in
S$. Therefore $\varphi$ induces a bijection on the set $\{v_s\,|\,
s\in S\}$. In other words $\varphi(\Delta^*)=\Delta^*$, and, since $\phi$ is
an isometry,
$\varphi(\Delta)=\Delta$ and the angle
between $\alpha_s$, $\alpha_r$ is preserved. That is, $\varphi$
induces an $u$-automorphism of the Coxeter graph of $W$, since the
order of $st$ in $W$ is entirely determined by the angle between
$\alpha_s$ and $\alpha_t$, and since $\kappa_s=\kappa_{t_s}$ for
all $s\in S$.
\end{proof}

\begin{rem} \textnormal{In the proof of the previous proposition, we have made use of the assumption
(stated when we introduced the root system $\Phi$) that all its roots
are of equal length.  Note that
this is not an important restriction, since any root
system can be rescaled to have all its roots of equal length.
For example, the vectors $(1,0)$ and $(-1,1)$ are simple roots for the
crystallographic root system $B_2$.  Instead of these, we would take
$(1,0)$ and $(-1/{\sqrt 2},1/\sqrt 2)$ as the simple roots. Attention, the assumption that
 the simple roots are of the same length does not imply that the fundamental weights are of the same length, see for instance in $A_3$}
\end{rem}

\begin{prop}\label{rem:automorphism_and_isometry}
  For every $u$-automorphism~$\mu$ of the Coxeter graph, there is a unique linear isometry~$\varphi_\mu$ that fixes~$P$ and $M(e)$ defined
   by $\varphi_\mu(\a_s):=\a_{\mu(s)}$ for every $s \in S$.
\end{prop}

\begin{proof}
The map $\varphi_\mu$
is well-defined since $\Delta$ is a basis of $V$. As~$\mu$ is an
   automorphism of the Coxeter graph and $\scalprod{\alpha_{\mu(s)}}{\alpha_{\mu(t)}}$ depends only on the order
   of~$st$, we have $\scalprod{\alpha_{\mu(s)}}{\alpha_{\mu(t)}}=\scalprod{\alpha_s}{\alpha_t}$ for $s,t\in S$.
   In other words~$\varphi_\mu$ is an isometry since $\Delta$ is a basis of~$V$.

From duality it is clear that  $v_s=\sum_{r\in
S}\scalprod{v_r}{v_s}\alpha_r$, for all $s\in S$. Also, the
matrices $[\scalprod{v_r}{v_s}]_{s,t}$ and
$[\scalprod{\a_r}{\a_s}]_{s,t}$ are inverse to each other and the
permutation $\mu\colon S\to S$ is such that
$[\scalprod{\a_{\mu(r)}}{\a_{\mu(s)}}]_{s,t}=[\scalprod{\a_r}{\a_s}]_{s,t}$.
Hence
$[\scalprod{v_{\mu(r)}}{v_{\mu(s)}}]_{s,t}=[\scalprod{v_r}{v_s}]_{s,t}$.
Thus, for $s\in S$ we have
$$
\varphi_\mu(v_s)=\sum_{r\in S}\scalprod{v_s}{v_r}\alpha_{\mu(r)}
=\sum_{r\in
S}\scalprod{v_{\mu(s)}}{v_{\mu(r)}}\alpha_{\mu(r)}=\sum_{r'\in S}
\scalprod{v_{\mu(s)}}{v_{r'}}\alpha_{r'}=v_{\mu(s)}.
$$
Now as $\kappa_s=\kappa_{\mu(s)}$ for all $s\in S$, $\varphi_\mu$
fixes $M(e)$, and therefore $P$.
\end{proof}

Similarly, for every isometry~$\varphi$ that fixes~$P$
and $M(e)$ there is an $u$-automorphism~$\mu$ such that
   $\varphi(v_s):=v_{\mu(s)}$ for every $s \in S$, by Proposition~\ref{prop:UnicityPerm}.

\begin{cor}\label{cor:group_action_and_isometry}
   Let~$\mu$ be an automorphism of the Coxeter graph of~$(W,S)$ and~$\varphi$ be a linear isometry that maps~$P$
   to itself and has~$M(e)$ as fixed point. Suppose that~$\mu$ and~$\varphi$ are related via $\varphi(v_s)=v_{\mu(s)}$
   for all $s \in S$. Then $\varphi=\varphi_\mu$ and $\mu$ is an
   $u$-automorphism. Moreover,
   \[
     \varphi\left(w(v_s)\right) = \left(\mu(w)\right)(v_{\mu(s)})
     \qquad\text{and}\qquad
     \varphi(\H_{(w,s)}) = \H_{(\mu(w),\mu(s))}
   \]
   for $w\in W$, $s\in S$. In particular, $\varphi(\Perm(W))=\Perm(W)$.
\end{cor}
\begin{proof}
   As $\Delta^*$ is a basis of $V$, $\varphi=\varphi_\mu$.
   Moreover, since $\varphi$ fixes $M(e)$ we have
   $$
     \sum_{\mu(s)\in S}\kappa_{\mu(s)}v_{\mu(s)}=\sum_{s\in
     S}\kappa_s v_s=M(e)=\varphi_\mu(M(e))= \sum_{s\in
     S}\kappa_sv_{\mu(s)}.
   $$
   By identification, we have $\kappa_s=\kappa_{\mu(s)}$ for all $s\in S$,  which proves that $\mu$ is an $u$-automorphism.

   We prove the first remaining claim by induction on the length of~$w$. If~$w=e$ the claim is~$\varphi(v_s)=v_{\mu(s)}$ and was shown
   in the proof of Proposition~\ref{prop:UnicityPerm}.
   Now assume $\ell(w)>0$. There is $t \in S$ such that $w=w^{\prime}t$ with $\ell(w^{\prime})<\ell(w)$. The action of~$t$
   on~$V$ is a reflection with reflection hyperplane~$H_t$. Hence we have
   \[
     t(v_s) = v_s - \frac{2\scalprod{v_s}{\a_t}}{\scalprod{\a_t}{\a_t}}\a_t.
   \]
   Also, since $\varphi$ is an isometry that maps~$\a_r$ to~$\a_{\mu(r)}$ and~$v_r$ to~$v_{\mu(r)}$, we have
   \begin{align*}
     \varphi_\mu(w(v_s)) &= \varphi_\mu\left(w^{\prime}\left(t(v_s)\right)\right) \\
                     &= \mu(w^{\prime})\left( v_{\mu(s)} - \frac{2\scalprod{v_{\mu(s)}}{\a_{\mu(t)}}}{\scalprod{\a_{\mu(t)}}{\a_{\mu(t)}}}\a_{\mu(t)}\right)\\
                     &= \mu(w^{\prime})\mu(t)\left(v_{\mu(s)}\right)
                      = \left(\mu(w^{\prime}t)\right)\left(v_{\mu(s)}\right)\\
                     &= \mu(w)\left(v_{\mu(s)}\right).
\end{align*}
We now prove the second claim.
\begin{align*}
    \varphi_\mu(\H_{(w,s)}) &= \{\varphi_\mu(v)\in V \,|\, \scalprod{v}{w(v_s)}\leq \scalprod{M(e)}{v_s}\}\\
                            &= \{v\in V \,|\, \scalprod{v}{\varphi_\mu(w(v_s))}\leq \scalprod{\varphi_\mu(M(e))}{\varphi_\mu(v_s)}\}\\
                            &= \{v\in V \,|\, \scalprod{v}{\mu(w)(v_{\mu(s)}))}\leq \scalprod{M(e)}{v_{\mu(s)}}\}\\
                            &= \H_{(\mu(w),\mu(s))}
\end{align*}
\end{proof}

\medskip
\begin{figure}[h]
  \begin{center}
  \begin{minipage}{0.95\linewidth}
     \begin{center}
     \begin{overpic}
        [width=10cm]{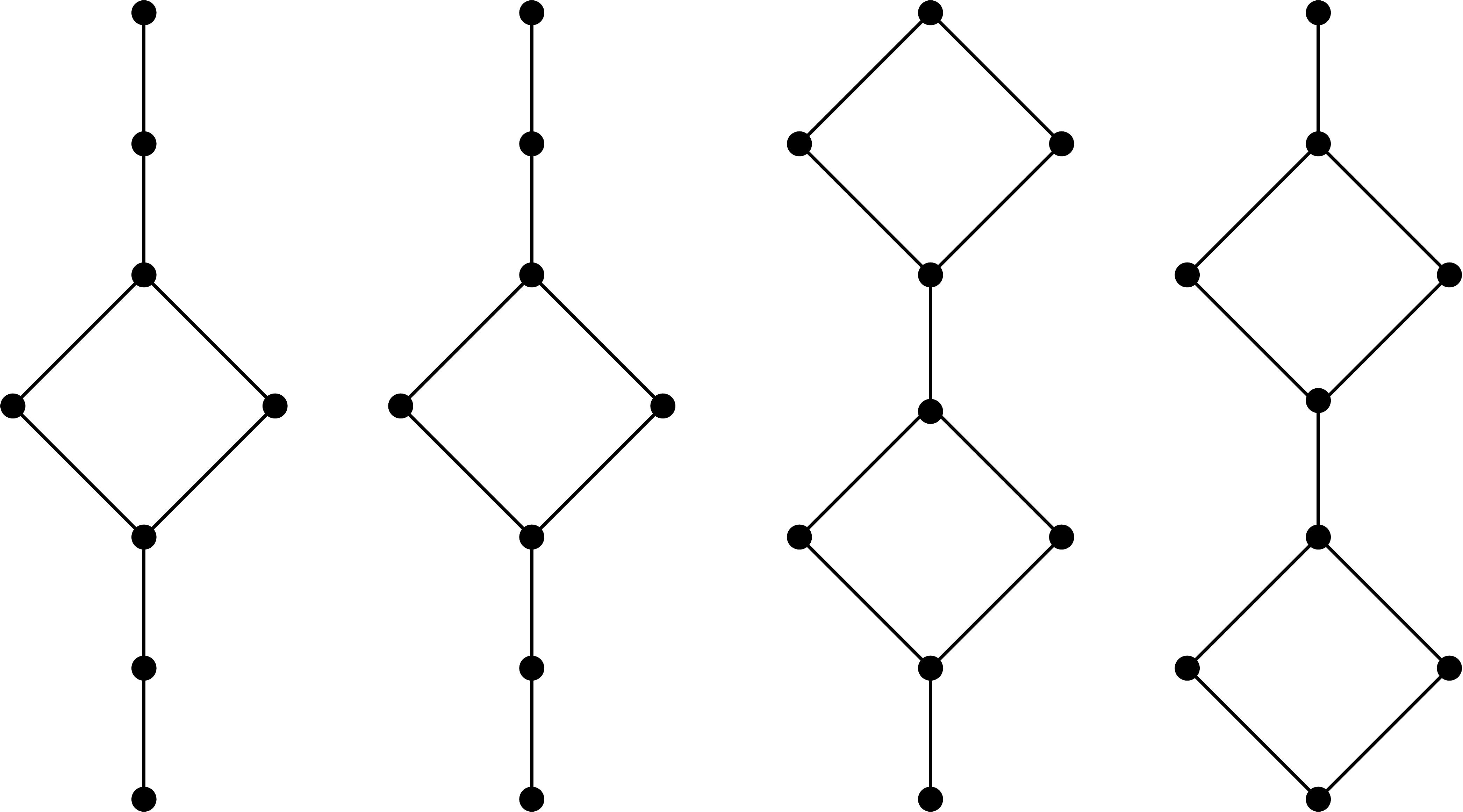}
        \put(11.5,0.5){\tiny $e$}
        \put(11.5,9){\tiny $s_1$}
        \put(11.5,18.5){\tiny $s_1s_2$}
        \put(2.5,27){\tiny $s_1s_2s_3$}
        \put(18,29){\tiny $s_1s_2s_1$}
        \put(11.5,36){\tiny $s_1s_2s_3s_1$}
        \put(11.5,45){\tiny $s_1s_2s_3s_1s_2$}
        \put(11.5,54.5){\tiny $s_1s_2s_3s_1s_2s_1$}
        \put(3,-4){$c=s_1s_2s_3$}
        \put(38,0.5){\tiny $e$}
        \put(38,9){\tiny $s_3$}
        \put(38,18.5){\tiny $s_3s_2$}
        \put(29,27){\tiny $s_3s_2s_1$}
        \put(47,27){\tiny $s_3s_2s_3$}
        \put(38,36){\tiny $s_3s_2s_1s_3$}
        \put(38,45){\tiny $s_3s_2s_1s_3s_2$}
        \put(38,54.5){\tiny $s_3s_2s_1s_3s_2s_3$}
        \put(30,-4){$c=s_3s_2s_1$}
        \put(65,0.5){\tiny $e$}
        \put(65,9){\tiny $s_2$}
        \put(56,18.5){\tiny $s_2s_3$}
        \put(74,18.5){\tiny $s_2s_1$}
        \put(65,27){\tiny $s_2s_1s_3$}
        \put(65,36){\tiny $s_2s_1s_3s_2$}
        \put(56,45){\tiny $s_2s_1s_3s_2s_1$}
        \put(72,47.5){\tiny $s_2s_1s_3s_2s_3$}
        \put(65,54.5){\tiny $s_2s_1s_3s_2s_1s_3$}
        \put(57,-4){$c=s_2s_1s_3$}
        \put(92,0.5){\tiny $e$}
        \put(83,9){\tiny $s_1$}
        \put(101,9){\tiny $s_3$}
        \put(92,18.5){\tiny $s_3s_1$}
        \put(92,27){\tiny $s_3s_1s_2$}
        \put(83,36){\tiny $s_3s_1s_2s_1$}
        \put(101,36){\tiny $s_3s_1s_2s_3$}
        \put(92,45){\tiny $s_3s_1s_2s_3s_1$}
        \put(92,54.5){\tiny $s_3s_1s_2s_3s_1s_2$}
        \put(84,-4){$c=s_3s_1s_2$}
     \end{overpic}
     \end{center}$ $
     \caption[]{There are four Coxeter elements in~$S_4$, each yields a distributive lattice~$\mathcal G_c$ of $c$-singletons.}
     \label{fig:dist_lattices}
  \end{minipage}
  \end{center}
\end{figure}
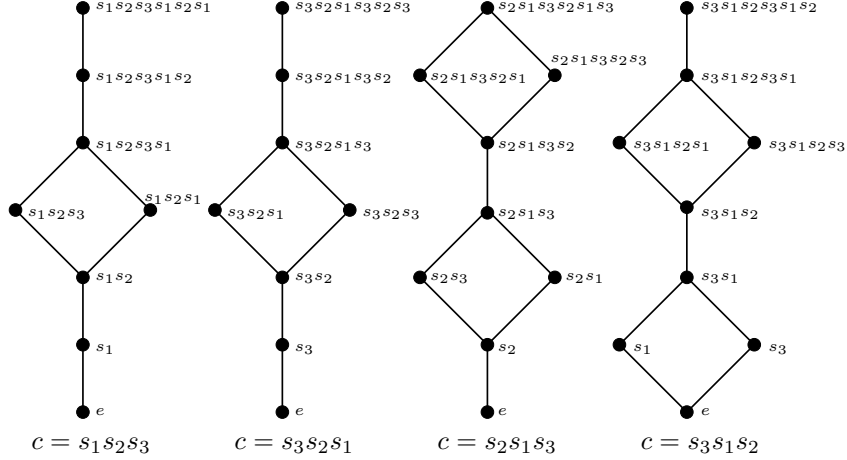

\medskip
The $c$-singletons form a distributive sublattice of the right
weak order, see~\cite{realisation2}. We denote the Hasse diagram
of this poset by~$\mathcal G_c$, see
Figure~\ref{fig:dist_lattices} for examples. They also form a
sublattice of the $c$-Cambrian lattice.

There is an important linear isometry on~$V$ that fixes~$P$ and
interchanges~$M(e)$ and~$M(w_0)$: The map~$g$ defined by $v
\longmapsto w_0(v)$.

Reading proves in \cite[Proposition~1.3]{reading3} that the map $w
\longmapsto ww_0$ is an anti-isomorphism from the $c$-Cambrian
lattice to the $c^{-1}$-Cambrian lattice. In particular, the map
$w \longmapsto ww_0$ is an anti-isomorphism between the lattices
of $c$-singletons and $c^{-1}$-singletons by restriction, that is,
$w$ is a $c$-singleton if and only if $ww_0$ is a
$c^{-1}$-singleton. Since the map $w\mapsto w_0ww_0$ is an
isomorphism from the $c$-Cambrian lattice to the
$w_0cw_0$-Cambrian lattice, the map $w\mapsto w_0w=(w_0ww_0) w_0$
is an anti-isomorphism from the $c$-Cambrian lattice to the
$w_0c^{-1}w_0$-Cambrian lattice. In other words, $\H_{(x,s)}$ is
$c$-admissible for all $s \in S$ if and only if $\H_{(w_0x,s)}$ is
$w_0c^{-1}w_0$-admissible for all $s \in S$. Therefore we obtain
the following proposition.

\begin{prop}\label{cor:wo} Let~$c$ be a Coxeter element of the finite Coxeter system~$(W,S)$.
 Then $\Ass_{c}(W)=g\left(\Ass_{w_0c^{-1}w_0}(W)\right)$, that is,
   the generalized associahedra~$\Ass_{c}(W)$ and~$\Ass_{w_0c^{-1}w_0}(W)$ are isometric.
\end{prop}

\noindent Let~$T$ be the reflections of~$W$ and~$I(w)$ be the
inversions of~$w\in W$ defined as
\[
  T:=\bigcup_{w\in W} wSw^{-1}
  \qquad\text{and}\qquad
  I(w):=\{t\in T\,|\, \ell(tw)<\ell(w)\}.
\]
A {\em parabolic subgroup} is a subgroup that is the conjugate of a standard parabolic subgroup of~$W$. Given the Coxeter system $(W,S)$ and a
parabolic subgroup $W'$, there is a
natural way to distinguish a set of {\it simple generators} of $W'$, see
\cite{reading2}.  We shall only make use of the case that
$W'$ is standard parabolic, in which case, the simple generators of
$W'$ are simply $W'\cap S$.

\begin{thm}\label{thm:CSingGraph}
   For $w\in W$ and any Coxeter element~$c$ of~$W$, the following statements are equivalent:
   \begin{compactenum}[(i)]
      \item $w$ and $ww_0$ are both $c$-singletons;
      \item $w \in \{ e,w_0 \}$;
      \item $wcw^{-1}$ is a Coxeter element of~$W$ and $w\mathcal G_c=\mathcal G_{wcw^{-1}}$;
      \item $w\mathcal G_c=\mathcal G_{c'}$ for some Coxeter element $c'$.
   \end{compactenum}
\end{thm}

\begin{proof} $ $
   $(i)\Rightarrow (ii)$:
   Suppose $w$ and $ww_0$ are $c$-singletons.
   A $c$-singleton~$u$ is $c$-sortable and
   $c$-antisortable by Proposition~2.2 of~\cite{realisation2}, that is, the element~$u$ is $c$-sortable and~$uw_0$ is
   $c^{-1}$-sortable. Hence~$w$ is $c$-sortable and $c^{-1}$-sortable. From \cite[Theorem~4.1]{reading2} we know
   that $g\in W$ is $c$-sortable and $c^{-1}$-sortable if and only if $I(g)\cap (W^{\prime}\setminus \{t_1\}) \neq \emptyset$
   implies $t_2\in I(g)$ for
any irreducible dihedral parabolic subgroup~$W^{\prime}$
   of~$W$ (that is $|W^{\prime}|>4$) with simple generators $t_1,t_2\in T$.

   Assume that $w \neq e$. There exists $s\in I(w)\cap S$. Pick $t\in S$ such that the standard parabolic
   subgroup~$W'$ generated by $\{s,t\}$ is dihedral and of cardinality~$>4$. We first show that $s,t \in I(w)$.
   We have to distinguish two cases:
   \begin{compactenum}[(1)]
      \item If $I(w)\cap (W'\setminus \{s\}) \neq \emptyset$ then $t\in I(w)$ because~$w$ is $c$-sortable
            and $c^{-1}$-sortable. Hence $s,t \in I(w)$.
      \item Assume $I(w)\cap (W'\setminus \{s\}) = \emptyset$. We first remark $I(ww_0)=I(w_0)\setminus I(w)$.
            Hence $I(ww_0)\cap (W'\setminus \{t\}) \neq \emptyset$ which implies $s \in I(ww_0)$ since $ww_0$
            is also $c$-sortable and $c^{-1}$-sortable. In particular, $s \in I(w) \cap I(ww_0)$ which is impossible.
   \end{compactenum}
   Since~$(W,S)$ is irreducible, the Coxeter graph associated to~$(W,S)$
   is connected. Now repeat this process along paths starting at~$s$ to conclude that $S\subseteq I(w)$. Hence~$w=w_0$.

   $(ii)\Rightarrow (iii)$:  For $w=e$ the result is clear. Recall that the conjugation by $w_0$ is an
   automorphism~$\varphi$ of the Coxeter system~$(W,S)$. So the $w_0cw_0$-factorization of $w_0$ is induced by~$\varphi$
   from the $c$-factorization of~$w_0$. The claim for $w=w_0$ follows.

   $(iii)\Rightarrow (iv)$: Set $c^{\prime} := wcw^{-1}$.

   $(iv)\Rightarrow (i)$: Since $e$ and~$w_0$ are $c^{\prime}$-singletons, we conclude that~$w^{-1}$ and~$w^{-1}w_0$ are
   both $c$-singletons. Apply $(i)\Rightarrow(ii)$ to deduce that $w^{-1}=e$ or $w^{-1}=w_0$. In particular,~$w$
   and~$ww_0$ are both $c$-singletons.
\end{proof}

\begin{rem}\label{rem:Cambrian}
   \textnormal{Two maximal cones~$C$ and~$C^{\prime}$ in the $c$-Cambrian fan~$\mathcal F_c$ are antipodal if~$C=-C'$.
       Theorem~\ref{thm:CSingGraph} implies that a pair of antipodal maximal cones that correspond to $c$-singletons is
       unique and the corresponding elements are~$e$ and~$w_0$.
      }
\end{rem}

\section{Proof of Theorem~\ref{thm:Main}}\label{se:Proof}

Assume there is an isometry~$\varphi$ on~$V$ such that $\Ass_{c_1}(W)=\varphi(\Ass_{c_2}(W))$.

\noindent
Let~$w$ be a $c_2$-singleton. Then
\begin{compactenum}
   \item $M(w) = \bigcap_{s\in S} H_{(w,s)}$ is a vertex of~$\Ass_{c_2}(W)$,
   \item $\varphi\left(M(w)\right) = \bigcap_{s\in S} \varphi(H_{(w,s)})$ is a vertex of~$\Ass_{c_1}(W)$,
   \item $\varphi\left( M(w)\right) = M(w^{\prime})$ for some $c_1$-singleton~$w^{\prime}$ since~$\varphi$ is an isometry.
\end{compactenum}
(For (3), note that $c$-singleton cones are the only cones in the Cambrian
fan which consist of a single chamber from the Coxeter fan, and thus an
isometry must take singleton cones to singleton cones.)

Apply these results
to~$w=e$ to obtain a $c_1$-singleton~$w^{\prime}_e$
with $M(w^{\prime}_e)=\varphi\left(M(e)\right)$. Moreover,
$w^{\prime}_e w_0$ is also a $c_1$-singleton with $M(w^{\prime}_e
w_0)=\varphi\left(M(w_0)\right)$. Hence $w^{\prime}_e \in \{
e,w_0\}$ by Theorem~\ref{thm:CSingGraph} and~$\varphi$ is a linear
isometry of~$V$ that fixes~$P$ and either fixes $M(e)$
and~$M(w_0)$ or interchanges~$M(e)$ and~$M(w_0)$. If~$\varphi$
fixes~$M(e)$ and~$P$ then there is an induced $u$-automorphism~$\mu$
of the Coxeter graph of~$(W,S)$ by
Proposition~\ref{prop:UnicityPerm} and $\mu(c_2)=c_1$. If
$\varphi$ interchanges~$M(e)$ and~$M(w_0)$ we consider
$\widetilde\varphi := g\circ\varphi$. We have
$\widetilde\varphi(\Ass_{c_2}(W)) = \Ass_{w_0c_1^{-1}w_0}(W)$ by
Proposition~\ref{cor:wo} and $\widetilde\varphi$ is an isometry
that fixes~$P$,~$M(e)$, and~$M(w_0)$. Hence $\widetilde\varphi$
induces an $u$-automorphism~$\mu$ of the Coxeter system~$(W,S)$ by
Proposition~\ref{prop:UnicityPerm} and $\mu(c_2)=w_0c_1^{-1}w_0$.

\medskip
Assume there is an $u$-automorphism~$\mu$ of the Coxeter
graph~$(W,S)$. Without loss of generality, we may assume
that~$\mu(c_2)=c_1$ because~$\Ass_{c}(W)$
and~$\Ass_{w_0c^{-1}w_0}(W)$ are isometric via~$g$ by
Proposition~\ref{cor:wo}.

We have to specify an isometry~$\varphi_\mu$ on~$V$ such that~$\varphi_\mu(\Ass_{c_2}(W))=\Ass_{c_1}(W)$.
This is done according to Proposition~\ref{rem:automorphism_and_isometry}: Define $\varphi_\mu:V\to V$ by
$\varphi_\mu(\a_s):=\a_{\mu(s)}$ for all~$s\in S$ or equivalently by
$\varphi_\mu(v_s):=v_{\mu(s)}$.

It remains to show that~$\varphi_\mu$ maps $c_2$-admissible
halfspaces to $c_1$-admissible halfspaces. From
Corollary~\ref{cor:group_action_and_isometry} we know how the
facet defining halfspaces~$\H_{(w,s)}$ are permuted by the
isometry~$\varphi_\mu$, that is, $\varphi_\mu(\H_{(w,s)}) =
\H_{(\mu(w),\mu(s))}$ for $w\in W$ and~$s\in S$. The
automorphism~$\mu$ on~$W$ preserves the length function~$\ell$, so
we have~$\mu(w_0)=w_0$ and any prefix of the $c$-factorization of
$w_0$ up to commutation is a prefix of the $\mu(c)$-factorization
of $w_0$ up to commutation. In other words, $\mu$ induces a
lattice isomorphism between the $c_2$-singletons and the
$\mu(c_2)$-singletons. Hence $\H_{(x,s)}$ is $c_2$-admissible if
and only if $\H_{(\mu(x),\mu(s))}$ is $\mu(c_2)$-admissible. This
shows that $\varphi_{\mu}(\Ass_c(W))=\Ass_{\mu(c)}(W)$ and ends
the proof of Theorem~\ref{thm:Main}.

\begin{proof}[Proof of Corollary~\ref{cor:Cambrian}] (2)
equivalent to (3) is Theorem~\ref{thm:Main}. (2) implies (1)
follows from the definition of normal fans.

Now, we show  (1) implies (3). As the $c$-Cambrian fan $\mathcal
F_c$ and the $c'$-Cambrian fan $\mathcal F_{c'}$ are isometric,
there is an isometry $\varphi$ such that the image under $\varphi$
of each cone in $\mathcal F_{c}$ is a cone of $\mathcal F_{c'}$.
By Remark~\ref{rem:Cambrian}, the pair of antipodal singleton cones
$C(e),C(w_0)$ corresponding to $e,w_0$ are the unique singleton antipodal
cones  in both Cambrian fans, and either $\varphi$ fixes them or
exchanges them. If $u$ is balanced, then apply
Corollary~\ref{cor:group_action_and_isometry} with $\mu$ to be the
conjugation by $w_0$ to obtain
\begin{equation}\label{eq:proofvalid}
 M(w_0)=\sum_{s\in S}\kappa w_0(v_s)=\sum_{s\in S}\kappa
 (-v_{w_0sw_0})=-M(e).
\end{equation}
As $C(e)$ is $\mathbb R_{>0}$-spanned by $\Delta^*$ and
$C(w_0)$ is $\mathbb R_{>0}$-spanned by $-\Delta^*$, either
$\varphi$ fixes $M(e)$ and $M(w_0)$ or interchanges them. In both
cases, $\varphi(P)=P$. So either $\varphi$ or $g\circ\varphi$
fixes $M(e)$ and $P$. Conclude by
Proposition~\ref{prop:UnicityPerm} as in the first part of the
proof of Theorem~\ref{thm:Main}.
\end{proof}

\section{The reducible case}

The reducible case does not follow immediately from an application of
the irreducible case.  Rather, one goes through the same steps as in the
proof of the irreducible case, but with some slight added complication.
We sketch the process below.

Let $\mathcal D$ denote the set of irreducible components of the
Coxeter graph of $(W,S)$.
For any $\mathcal A\subset \mathcal D$, let $w_{\mathcal A}$ be the
longest word for the subgroup generated by the components in $\mathcal A$.
Let $L=\{ w_{\mathcal A}\mid \mathcal A \subset \mathcal D\}$.   All
the elements of $L$ are $c$-singletons (for any $c$).

Up to just before Proposition~\ref{cor:wo}, the argument goes in exactly the same way.
Then, instead of
constructing a single isometry $g$, we construct an isometry
$g_{\mathcal A}$ for each $\mathcal A\subset \mathcal D$, defining
$g_{\mathcal A}(v)=w_{\mathcal A}(v)$.
Write $c^{\mathcal A}$ for the Coxeter element obtained from $c$ by
reversing the order of the reflections in $c$ coming from components in
$\mathcal A$.  The generalization of Proposition~\ref{cor:wo} then asserts that
$\Ass_c(W)=g_{\mathcal A}(\Ass_{w_{\mathcal A}c^{\mathcal A}w_{\mathcal A}}(W))$ and, in
particular, these associahedra are isometric.
Theorem~\ref{thm:CSingGraph} goes through with condition (ii) replaced by the condition that
$w\in L$.


\end{document}